\documentclass[11pt,reqno]{amsart}
\usepackage{yhmath}
\usepackage{fancyhdr}
\usepackage{times}
\usepackage[T1]{fontenc}
\usepackage{mathrsfs}
\usepackage{latexsym}
\usepackage[dvips]{graphics}
\usepackage{epsfig}
\usepackage{hyperref, flafter}
\usepackage{amsmath,amsfonts,amsthm,amssymb,amscd}
\usepackage{color}
\usepackage{geometry}               
\usepackage{epstopdf}
\usepackage{verbatim}

\allowdisplaybreaks

\newtheorem{thm}{Theorem}[section]

\newtheorem*{thm*}{Theorem}
\newtheorem*{con*}{Conjecture}
\newtheorem*{lem*}{Lemma}
\newtheorem{prop}[thm]{Proposition}

\theoremstyle{definition}
\newtheorem{rmk}{Remark}

\def\t{\theta}

\begin{document}
	
	\title[log derivative CUE]
	{On the logarithmic derivative of characteristic polynomials for random unitary matrices}
	
	\author{Fan Ge}

	\address{Department of Mathematics, William \& Mary, Williamsburg, VA, United States}
	
	\email{ge@wm.edu}
	
	\baselineskip=15pt
	
	\begin{abstract}
		Let $U\in U(N)$ be a random unitary matrix of size $N$, distributed with respect to the Haar measure on $U(N)$.  Let $P(z)=P_U(z)$ be the characteristic polynomial of $U$.  We prove that for $z$ close to the unit circle, $ \frac{P'}{P}(z) $ can be approximated using zeros of $P$ very close to $z$, with a typically controllable error term. This is an analogue of a result of Selberg for the Riemann zeta-function. We also prove a mesoscopic central limit theorem for $ \frac{P'}{P}(z) $ away from the unit circle, and this is an analogue of a result of Lester for zeta.
	\end{abstract}
	
	\maketitle
	
	\section{Introduction}

	Let $\zeta(s)$ be the Riemann zeta-function, and let $\rho=\beta+i\gamma$ denote a generic nontrivial zero of zeta. 
	A beautiful result of Selberg~\cite{selbergzeta'} says that for $s=\frac{1}{2}+it$ with $t\in [T, 2T]$ and $s\ne \rho$ we can write
	\begin{align}\label{eq sel}
		\frac{\zeta'}{\zeta}(s) = \sum_{|t-\gamma|<\frac{1}{\log T}} \frac{1}{s-\rho} + D,
	\end{align}
	where the error term $D$ is in terms of an explicit Dirichlet polynomial and satisfies
	\begin{align*}
		\frac{1}{T}	\int_{T}^{2T} |D|^{2K} dt \ll_K \log^{2K} T
	\end{align*}
	for all positive integers $K$. In other words, $\zeta'/\zeta(s)$ can be approximated using zeros very close to $s$, with a typically controllable error. Radziwi{\l}{\l}~\cite{radziwill2014gaps} observed that Selberg's argument with some modification also gives that for every constant $0<c\le 1$ and $s$ close to the critical line ($0\le \Re(s)-1/2\ll 1/\log T$, say),
	\begin{align}\label{eq rad}
		\frac{\zeta'}{\zeta}(s) = \sum_{|t-\gamma|<\frac{c}{\log T}} \frac{1}{s-\rho} + E,
	\end{align}
	with
	\begin{align*}
		\frac{1}{T}	\int_{T}^{2T} |E|^{2K} dt \ll_K \left(\frac{\log T}{c}\right)^{2K}.
	\end{align*}
	(Radziwi{\l}{\l}'s paper assumed the Riemann Hypothesis mainly for other purposes; for \eqref{eq rad} alone one can show it holds unconditionally.) For applications of these results, see for example Selberg~\cite{selbergzeta'}, Radziwi{\l}{\l}~\cite{radziwill2014gaps}, and Ge~\cite{ge2023zeta}.
	
	It is well known that characteristic polynomials of the circular unitary ensemble (CUE) models the Riemann zeta-function $\zeta(s)$, with the matrix size $N$ about the same as $\log T$.  Our first result in this paper is a CUE analogue of~\eqref{eq sel} and~\eqref{eq rad}. Throughout, let  $U(N)$ be the set of unitary matrices of size $N$, equipped with Haar measure. For $U\in U(N)$, write 
	\begin{align*}
		P(z)=P_U(z)=\prod_{j=1}^{N} (z-z_j)
	\end{align*}
	for the characteristic polynomial of $U$, where 
	\begin{align*}
		z_j=e^{i\t_j} \text{ with } -\pi<\t_j\le \pi, \text{ for } j=1, 2,..., N.
	\end{align*}
	
	\begin{thm}\label{thm sel cue}
		Let $0<c\le 1$  be a constant. For $1-\frac{1}{N} \le z \le 1$ and $z\ne z_j$ we have
		\begin{align*}
			\frac{P'}{P}(z) = \sum_{|\t_j|<\frac{c}{N}} \frac{1}{z-z_j} + \mathcal E,
		\end{align*}
		where the error term $\mathcal E= \sum_{|\t_j| \ge \frac{c}{N}} \frac{1}{z-z_j}$ satisfies
		\begin{align*}
			\mathbb E |\mathcal E|^{2K} \ll_K \left(\frac{N}{c}\right)^{2K}
		\end{align*}
		for every positive integer $K$. Here the expection $\mathbb E$ is over $U(N)$ with respect to the Haar measure.
	\end{thm}

\begin{rmk}
Since the Haar density of a configuration of eigenvalues in CUE is invariant under rotations, Theorem~\ref{thm sel cue} holds not just for real $z$ but for all $1-\frac{1}{N} \le |z| \le 1$ with obvious modifications in the statement.
\end{rmk}

	\begin{rmk}
		Bailey, Bettin, Blower, Conrey, Prokhorov, Rubinstein and Snaith~\cite{bailey2019mixedmoments}  proved that for $K\in\mathbb Z^+$
		\begin{align}\label{thm Bailey et al}
			\int_{U(N)} \left|\frac{P'}{P}\Big(1-\frac{a}{N}\Big)\right|^{2K} dU \sim \binom{2K-2}{K-1} \frac{N^{2K}}{(2a)^{2K-1}}
		\end{align}
		where $a\to 0$ as $N\to \infty$. They also conjectured a similar asymptotic for zeta, and this was studied in~\cite{ge2023zeta}.
In \cite{alvarez2020moments} Alvarez and Snaith proved analogous results of~\eqref{thm Bailey et al} for orthogonal and symplectic random matrices. More recently, Alvarez, Bousseyroux and Snaith~\cite{alvarez2023noninteger} extended the corresponding result for the odd orthogonal ensemble to non-integer moments (namely, $K\in\mathbb R^+$). In~\cite{ge2023real} we shall apply our Theorem~\ref{thm sel cue} to obtain asymptotics for real moments analogues of~\eqref{thm Bailey et al} in unitary, even orthogonal, and symplectic ensembles.
	\end{rmk}

	We also  investigate the value distribution of $ 	\frac{P'}{P}(z)  $ in the mesoscopic range away from the unit circle, and prove that it obeys a central limit theorem (CLT). To put it in context,  Selberg's central limit theorem states that, roughly speaking, the value distribution of the vector
	\begin{align*}
		\frac{1}{\sqrt{\frac{1}{2}\log\log T}}	\left( \Re \log \zeta (1/2+it) ,  \Im \log \zeta (1/2+it) \right), \quad t\in [T, 2T]
	\end{align*}
	converges to a  normal random vector $(X, Y)$ as $T$ tends to infinity. Here $X$ and $Y$ are independent and both have mean $0$ and variance $1$, and $t$ is drawn uniformly from $[T, 2T]$. The result also holds for $\sigma+it$ in place of $1/2+it$ if $\sigma$ is close to $1/2$, but the variance changes accordingly. See Tsang's thesis~\cite{tsang1984thesis} for a detailed description and proof.  The imaginary part of $\log \zeta(1/2+it)$ is related to the counting of zeta zeros via the Riemann-von Mangoldt formula; in this direction Selberg's CLT for $ \Im \log \zeta (1/2+it)$ is related to macroscopic and mesoscopic CLTs for value distributions (again $t$ being drawn uniformly in $[T, 2T]$) of quantities of the form $\sum_\rho \eta(\Delta\cdot (\rho-t))$, where the sum is over zeta zeros, $\eta$ is a suitable function, and $\Delta$ is a scaling factor. Here macroscopic scale refers to the case when $\Delta$ is of constant size, and mesoscopic scale is the case when $\Delta\to \infty$ with $T$ but $\Delta=o(\log T)$. See Fujii~\cite{fuji1974zeros}, Bourgade and Kuan~\cite{bourgade2014clt}, Rodgers~\cite{rodgers2014clt}, and Maples and Rodgers~\cite{rodgers2015clt2}. Another related result is a mesoscopic CLT of Lester~\cite{lester2014distribution} for the logarithmic derivative of zeta away from the critical line, who proved that when $t$ ranges from $T$ to $2T$  the real part and the imaginary part of
	\begin{align*}
		\frac{\zeta'}{\zeta} \left(\frac{1}{2} + \frac{\psi(T)}{\log T} + it\right)
	\end{align*}
	are close to independent normal with mean $0$ and variance 
	\begin{align*}
		V_{zeta}=	\frac{1}{2} \sum_{n=2}^{\infty} \frac{\Lambda(n)^2}{n^{1+\frac{2 \psi(T)}{\log T}}},
	\end{align*}
	provided that $\psi(T)=o(\log T)$ and $\psi(T)\to \infty$ with $T$.
	A standard calculation shows that 
	\begin{align}\label{eq lester v}
		V_{zeta} \sim	\frac{1}{8} \left(\frac{\log T}{\psi(T)}\right)^2.
	\end{align}
	Analogues of Selberg's CLT as well as macroscopic and mesoscopic CLTs for sums over zeta zeros are also known in the CUE setting; see Keating and Snaith~\cite{keating2000momentszeta}, Hughes, Keating and O'Connell~\cite{hughes2001clt}, Bourgade~\cite{bourgade2010cltaway}, Szeg\"o~\cite{szego1952stronglimit}, Wieand~\cite{wieand2002},  Diaconis and Evans~\cite{diaconis2001clt}, and Soshnikov~\cite{soshnikov2000combinatorial}. Related to Lester's result, in the CUE setting there is a mesoscopic CLT away from the real line  proved by Chhaibi, Najnudel and Nikeghbali~\cite{chhaibi2017ratios} for a limiting object of the characteristic polynomials. Our next result is a CUE analogue of Lester's CLT.

	\begin{thm}\label{thm clt}
	Let $L=L_N$ be a quantity such that $L_N=o(N)$ and $L_N\to \infty$ as $N\to \infty$. Then with $U$ drawn from $U(N)$ with respect to the Haar measure, the random vector
	\begin{align*}
			\left(\Re \left[ \frac{L}{N} \cdot 	\frac{P'}{P}\left(1-\frac{L}{N}\right) \right], \Im \left[ \frac{L}{N} \cdot 	\frac{P'}{P}\left(1-\frac{L}{N}\right) \right] \right)
	\end{align*}
		converges in distribution to a normal vector $(G, H)$ as $N\to\infty$, where $G$ and $H$ are independent normals with mean $0$ and variance $1/8$.
	\end{thm}
	Therefore,  for large $N$ the real part and the imaginary part of $\frac{P'}{P}\left(1-\frac{L}{N}\right)$ are close to independent normal random variables with mean $0$ and variance $V$ with
	\begin{align}
		V=\frac{1}{8} \left(\frac{N}{L}\right)^2.
	\end{align}
	This variance agrees  with Lester's variance~\eqref{eq lester v} in the zeta case.

	\section{Proof of Theorem~\ref{thm sel cue}}
	
	Define
	\begin{align}\label{eq z0}
		z_0 = 1 - \frac{1}{N}
	\end{align}
	and let
	\begin{align}\label{eq z}
		z_0 \le z\le 1.
	\end{align}
	We write 
	\begin{align*}
		\frac{P'}{P}(z) = \sum_{|\t_j|<\frac{c}{N}} \frac{1}{z-z_j} + X_1 +X_2 -X_3,
	\end{align*}
	where
	\begin{align*}
		X_1 &= 	\frac{P'}{P}(z_0) = \sum_{j=1}^N \frac{1}{z_0-z_j}, \\
		X_2 &=  \sum_{|\t_j|\ge \frac{c}{N}} \left( \frac{1}{z-z_j} - \frac{1}{z_0-z_j} \right),\\
		X_3 &= \sum_{|\t_j|<\frac{c}{N}} \frac{1}{z_0-z_j}.
	\end{align*}
	Theorem~\ref{thm sel cue} follows from the next three propositions.
	\begin{prop} \label{prop 1}
	For $K\in \mathbb Z^+$ we have	$	\mathbb E |X_1|^{2K} \ll _K N^{2K}$.
	\end{prop}
	
	\begin{prop} \label{prop 2}
		For $K\in \mathbb Z^+$ we have		$	\mathbb E |X_2|^{2K} \ll _K \left(\frac{N}{c}\right)^{2K}$.
	\end{prop}
	
	\begin{prop} \label{prop 3}
		For $K\in \mathbb Z^+$ we have		$	\mathbb E |X_3|^{2K} \ll _K N^{2K}$.
	\end{prop}
	
	\subsection{Proof of Proposition~\ref{prop 1}}
	
	We will need the following ratios formula from Conrey and Snaith~\cite{conrey2008correlations}.
	\begin{thm}\label{thm ratio} (Conrey and Snaith.)
		If $\Re \alpha_j>0$ and $ \Re \beta_j
		>0$ for  $\alpha_j\in A$ and $\beta_j\in B$, then $J(A;B)=J^*(A;B)$  where
		\begin{eqnarray*}
			J(A;B)&:=& \int_{U(N)}\prod_{\alpha\in A}
			(-e^{-\alpha})\frac{P_U'}{P_U}(e^{-\alpha})\prod_{\beta\in
				B} (-e^{-\beta})\frac{P_{U^*}'}{P_{U^*}}(e^{-\beta})~ dU
			,
		\end{eqnarray*}
		\begin{eqnarray*} &&J^*(A;B):= \nonumber \\
			&&\qquad\qquad\sum_{S\subset A,T\subset B\atop
				|S|=|T|}e^{-N(\sum_{\hat \alpha\in S} \hat \alpha
				+\sum_{\hat{\beta}\in T}\hat\beta)} \frac{Z(S,T)Z(S^-,T^-)} {
				Z^{\dagger}(S,S^-)Z^{\dagger}(T,T^-)} \sum_{{(A-S)+ (B-T)\atop =
					U_1+\dots + U_R}\atop |U_r|\le 2}\prod_{r=1}^R H_{S,T}(U_r),
		\end{eqnarray*}
		and
		\begin{equation*}\label{eqn:Hrmt}
			H_{S,T}(W)=\left\{\begin{array}{ll} \sum_{\hat \alpha\in
					S}\frac{z'}{z}(\alpha-\hat{\alpha})-\sum_{\hat\beta\in T}
				\frac{z'}{z}(\alpha +\hat \beta) &\mbox{ if $W=\{\alpha\}\subset
					A-S$}
				\\
				\sum_{\hat\beta\in T}\frac{z'}{z}(\beta-\hat
				\beta)-\sum_{\hat\alpha\in S} \frac{z'}{z}
				(\beta+\hat\alpha) &\mbox{ if  $W=\{\beta\}\subset B-T$}\\
				\left(\frac{z'}{z}\right)'(\alpha+\beta) & \mbox{ if
					$W=\{\alpha,\beta\}$ with $
					{\alpha \in A-S, \atop \beta\in B-T}$}\\
				0&\mbox{ otherwise}.
			\end{array}\right.
		\end{equation*}
		Here $z(x)=(1-e^{-x})^{-1}$, $S^-=\{-s:s\in S\}$ (similarly for $T^-$) and  $Z(A,B)=\prod_{\alpha\in A\atop\beta\in B}z(\alpha+\beta)$,
		with the dagger on $Z^\dagger(S,S^-)$ imposing the additional
		restriction that a factor $z(x)$ is omitted if its argument is
		zero.
	\end{thm}
	
	Here the $U^*$ is the conjugate transpose of $U$, and it is easy to see that
	\begin{align*}
		\mathbb E \left|\frac{P_U'}{P_U}(z)\right|^{2K} = \int_{U(N)} \left(\frac{P_U'}{P_U}(z)\right)^K \left( \frac{P_{U^*}'}{P_{U^*}}(\overline z) \right)^K ~ dU.
	\end{align*}
	One may try to apply the above theorem by letting all the $\alpha$'s in $A$ equal and similar for $\beta$'s in $B$. However, this will cause complications because many terms in the $J^*$ will have poles. To get around  this, we will apply the ratios formula to a discretized set of points, as follows.
	Let $r=\frac{1}{2N}$,  $r_k=\frac{r}{2^k}$, and define 
	$$G(x_1, ..., x_K) = \int_{U(N)} \frac{P_U'}{P_U}(x_1)  \frac{P_{U^*}'}{P_{U^*}}(x_1)\cdots \frac{P_U'}{P_U}(x_K)  \frac{P_{U^*}'}{P_{U^*}}(x_K) ~ dU. $$
	Since $z_0$ is real, to prove the proposition it suffices to show that $G(z_0, ..., z_0)\ll_K N^{2K}.$
	Observe that $G$ is holomorphic in each variable in a certain domain around $z_0$, so applying the maximum modulus principle to $G(x, z_0, ..., z_0)$ as a function of $x$, we see that there exists a $w_1$ on the circle $C(z_0, r_1)$ centered at $z_0$ with radius $r_1$ such that  $$|G(w_1, z_0, ..., z_0)| \ge |G(z_0, z_0, ..., z_0)| .$$  Next,  apply the maximum modulus principle to $G(w_1, x, z_0, ..., z_0)$ as a function of $x$, and we see that there exists a $w_2$ on the circle $C(z_0, r_2)$  such that $$ |G(w_1, w_2, z_0, ..., z_0)| \ge |G(w_1, z_0, ..., z_0)|  .$$  Repeat this and we conclude that there are points $w_1,...,w_k$ with $w_i$ on the circle $C(z_0, r_i)$ for each $i$, such that $$ |G(w_1,...,w_k)| \ge |G(z_0,..., z_0)|. $$
	By the definition~\eqref{eq z0} of $z_0$, it is easy to see that $\min_{i,j} |w_i-w_j| \gg_K r$. Since the function $\log(1-x)$ is close to $-x$ around $x=0$, we conclude that the points $\log(w_k)$'s are also well-spaced by a distance $\gg r \gg 1/N$. 
	We can now apply the ratios formula to the $w_k$'s. In the formula for $J^*$, from the above observation that these points are well-spaced, it is easy to see that all the $Z(S,T)$ and $Z(S^-, T^-)$ factors on the numerator are each bounded by $N^{|S|^2}$, and all the 	$Z^{\dagger}(S,S^-)$ and $Z^{\dagger}(T,T^-)$ factors on the denominator are each bounded from below by $N^{|S|^2-|S|}$. (Indeed, for $x$ close to $0$ the function $z(x)$ is close to $1/x$, and thus, the $z$ factors in the $Z$ and $Z^{\dagger}$ functions are $\gg$ and $\ll N$ by our choice of the $w_k$'s. The claimed bounds for $Z$ and $Z^{\dagger}$ follow by noticing that there are $|S|^2$ factors in $Z(S,T)$ and $|S|^2 - |S|$ factors in $Z^{\dagger}(S, S^-)$.) Similarly, all the $H_{S,T}$ factors are bounded by $N^{2K-2|S|}$. Collecting these estimates we conclude that $J^* \ll _K N^{2K}$, and this finishes the proof of Proposition~\ref{prop 1}.

	\subsection{Proof of Proposition~\ref{prop 2}}
Recall that
$
z_0 = 1 - \frac{1}{N}
$
and
$
z_0 \le z\le 1.
$
Observe that
	\begin{align*}
		|z_0-z_j| \le |z_0-z| + |z-z_j| \ll 1/N + |z-z_j|
	\end{align*}
	and that for $|\t_j| \ge c/N$
	\begin{align*}
		|z-z_j| \gg c/N.
	\end{align*}
It follows that
\begin{align*}
	c	|z_0-z_j| \ll c/N + c |z-z_j| \ll 	|z-z_j| + c |z-z_j| ,
\end{align*}
 and thus, for $0<c\le 1$  we have
	\begin{align*}
	c	|z_0-z_j| \ll  |z-z_j|.
\end{align*}
Using this we see that
	\begin{align*}
		X_2 &=  \sum_{|\t_j|\ge \frac{c}{N}} \left( \frac{1}{z-z_j} - \frac{1}{z_0-z_j} \right) \\
		& =  \sum_{|\t_j|\ge \frac{c}{N}}  \frac{z_0-z}{(z-z_j)(z_0-z_j)} \\
		& \ll \frac{1}{cN} \sum_{|\t_j|\ge \frac{c}{N}}  \frac{1}{|z_0-z_j|^2} \\
		& \ll  \frac{1}{cN} \sum_{j=1}^N  \frac{1}{|z_0-z_j|^2} ,
	\end{align*}

	We claim that
	\begin{align} \label{eq claim}
		\sum_{j=1}^N  \frac{1}{|z_0-z_j|^2} \ll N(N+|X_1|).
	\end{align}
	From this together with Proposition~\ref{prop 1} we will obtain Proposition~\ref{prop 2}. 
	
	To prove the claim~\eqref{eq claim}, we would like to connect its left-hand side with 
	\begin{align*}
		X_1 = 	\frac{P'}{P}(z_0) = \sum_{j=1}^N \frac{1}{z_0-z_j},
	\end{align*}
	and a natural way is to take the real or imaginary part of $X_1$. However, unlike in the zeta case, after taking real or imaginary parts the summands may have different signs or even be $0$, and this prevents us from controlling the size of the left-hand side of~\eqref{eq claim}. To get around this,
	we  introduce a change of variable 
	\begin{align*}
		P(z)=P(e^s)=Q(s)=\prod_{j=1}^{N} (e^s - e^{i\t_j}).
	\end{align*}
	We then consider a variant of the Weierstrass-Hadamard factorization for each $e^s - e^{i\t_j}$ and first write for purely imaginary $s$ that
	\begin{align*}
		e^s - e^{i\t_j} & = e^{\frac{s+i\t_j}{2}}\cdot (2i)\cdot \sin \left(\frac{-is-\t_j}{2}\right) \\
		& = e^{\frac{s+i\t_j}{2}}\cdot (s-i\t_j)\cdot \prod_{n\ne 0} \left(1+\frac{is+\t_j}{2n\pi}\right).
	\end{align*}
	It follows that
	\begin{align}\label{eq hadamard analogue}
		\notag	\frac{Q'}{Q}(s) & = \sum_{j=1}^N \frac{(	e^s - e^{i\t_j})'}{	e^s - e^{i\t_j}} \\
		\notag	& =  \sum_{j=1}^N  \left( \frac{1}{2} + \frac{1}{s-i\t_j} +\sum_{n\ne 0} \frac{1}{-i2n\pi + s - i\t_j}\right) \\
		\notag	& =  \sum_{j=1}^N  \left( \frac{1}{2} +\sum_{n\in \mathbb Z} \frac{1}{-i2n\pi + s - i\t_j}\right) \\
		& = \frac{N}{2} + \sum_{j=1}^N \sum_{n\in \mathbb Z} \frac{1}{s-i(\t_j+2n\pi)}.
	\end{align}
	Remark: One may view this formula as an analogue of the Hadamard fraction formula for zeta, where the $N/2$ in~\eqref{eq hadamard analogue} corresponds to the contribution of trivial zeta zeros which is about $-\frac{\log T}{2}$. The difference of the sign comes from the fact that $s$ is to the left of the `critical line' in the CUE case.

	These equations extend to all $s\in \mathbb C$ (except at poles), and thus,
	\begin{align*}
		\Re	\left(\frac{Q'}{Q}(s)\right) & = \frac{N}{2} + \sum_{j=1}^N \sum_{n\in \mathbb Z} \frac{\Re (s)}{|s-i(\t_j+2n\pi)|^2}.
	\end{align*}
	Set $s=s_0$ in the above equation, where $e^{s_0}=z_0$, so that $s_0=\log(1-\frac{1}{N})$ is about $-1/N$. We have
	\begin{align*}
		0<  \sum_{j=1}^N  \frac{-s_0}{|s_0-i\t_j|^2} & \le \sum_{j=1}^N  \sum_{n\in\mathbb Z} \frac{-s_0}{|s_0-i(\t_j+2n\pi)|^2} \\
		& = \frac{N}{2} - \Re	\left(\frac{Q'}{Q}(s_0)\right) \\
		& \le \frac{N}{2} + 	\left|\frac{Q'}{Q}(s_0)\right|.
	\end{align*}
	Now observe that 
	\begin{align*}
		|z_0-z_j| \gg |s_0 - i\t_j|
	\end{align*}
	and
	\begin{align*}
		\left|	\frac{Q'}{Q}(s_0) \right| = \left|\frac{P'(e^{s_0})\cdot e^{s_0}}{P(e^{s_0})} \right| \ll \left|\frac{P'(z_0)}{P(z_0)} \right| = |X_1| .
	\end{align*}
	Combine these with the above and we see that
	\begin{align*}
		\frac{1}{N} \sum_{j=1}^N  \frac{1}{|z_0-z_j|^2}  
		& \ll 	\frac{1}{N} \sum_{j=1}^N  \frac{1}{|s_0-i\t_j|^2}  \\
		& \ll \sum_{j=1}^N  \frac{-s_0}{|s_0-i\t_j|^2} \\
		& \ll N +	\left|\frac{Q'}{Q}(s_0)\right| \\
		& \ll N +|X_1|,
	\end{align*}
	proving~\eqref{eq claim}.
	
	\subsection{Proof of Proposition~\ref{prop 3}}
	
	We have 
	\begin{align*}
		X_3 & = \sum_{|\t_j|<\frac{c}{N}} \frac{1}{z_0-z_j}  \ll \left(\sum_{|\t_j|<\frac{c}{N}} 1 \right) \cdot N \le \frac{1}{N} \cdot \left(\sum_{|\t_j|<\frac{1}{N}} 1 \right) \cdot N^2.
	\end{align*}
	To bound the number of $|\t_j|<1/N$, we observe that
	\begin{align*}
		\sum_{j=1}^N  \frac{1}{|z_0-z_j|^2} \gg \sum_{|\t_j|<1/N}  \frac{1}{|z_0-z_j|^2} \gg N^2 \cdot \left(\sum_{|\t_j|<\frac{1}{N}} 1 \right).
	\end{align*}
	It follows that
	\begin{align*}
		X_3   \ll \frac{1}{N} \cdot \sum_{j=1}^N  \frac{1}{|z_0-z_j|^2}  \ll N+|X_1|,
	\end{align*}
	where the last inequality is by~\eqref{eq claim}.  Proposition~\ref{prop 3} now follows from Proposition~\ref{prop 1}.
	

	\section{Proof of Theorem~\ref{thm clt}}
	
	In view of Lester's result, we expect that the variance of the complex random variable 
	\begin{align*}
		\frac{P'}{P}\left(1-\frac{L}{N}\right) = \sum_{j=1}^{N} \frac{1}{(1-L/N)-e^{i\t_j}}.
	\end{align*}
	is a constant times $(N/L)^2$. Therefore, we introduce the following rescaling
	\begin{align*}
		&	f(\t)=f_{N,L}(\t) =\frac{L}{N} \cdot \frac{1}{(1-L/N)-e^{i\t}}, \\
		&	g(\t)=g_{N,L}(\t)= \Re f(\t), \\
		&		h(\t)=h_{N,L}(\t)= \Im f(\t),
	\end{align*}
	and denote
	\begin{align*}
		&	S_N(f)=\sum_{j=1}^{N} f(\t_j) =	 \frac{L}{N} \cdot 	\frac{P'}{P}\left(1-\frac{L}{N}\right), \\
		&	S_N(g)=\sum_{j=1}^{N} g(\t_j) = \Re \left[ \frac{L}{N} \cdot 	\frac{P'}{P}\left(1-\frac{L}{N}\right) \right], \\
		&	S_N(h)=\sum_{j=1}^{N} h(\t_j) = \Im \left[ \frac{L}{N} \cdot 	\frac{P'}{P}\left(1-\frac{L}{N}\right)\right].
	\end{align*}

	We compute the characteristic function (ch.f.) 
	\begin{align*}
		\phi_N(u,v)= \mathbb E e^{i(uS_N(g)+vS_N(h))}
	\end{align*}
	of the random vector $(S_N(g), S_N(h))$. We will show that 
	\begin{align}\label{eq ch.f. convg}
		\phi_N(u,v) \to e^{-\frac{1}{2}\cdot \frac{u^2+v^2}{8}} \quad \text{ as } N\to \infty,
	\end{align}
	for every $(u,v)\in \mathbb R^2$, and this will prove Theorem~\ref{thm clt} according to the convergence theorem for random vectors (see for example Theorem 3.10.5 in Durrett~\cite{durrett2019book}).
	To prove~\eqref{eq ch.f. convg}, it suffices to show that for every $(u,v)\in \mathbb R^2$,  the real  random variable 
	\begin{align}\label{eq univariate}
		uS_N(g)+vS_N(h) \longrightarrow Normal\left(0, \frac{u^2+v^2}{8}\right) \text{ in distribution,}
	\end{align}  
	for this will imply pointwise convergence of the ch.f. of $ uS_N(g)+vS_N(h) $, thus in particular its ch.f. evaluated at $1$, which gives~\eqref{eq ch.f. convg}.
	
	The main tool we use to prove \eqref{eq univariate} is the following result of Soshnikov~\cite{soshnikov2000combinatorial}, which is a combination of Lemma 1 and the main combinatorial lemma in that paper. 
	\begin{prop} \label{prop sosh}
		(Soshinikov.)
		Let $F(\t)$ be a real-valued function on the unit circle with continuous derivative and satisfy 
		$$\sum_{k\in\mathbb Z} |k||\hat F(k)|^2 < \infty,$$ 
		where $$\hat{F}(k)=\frac{1}{2\pi} \int_{0}^{2\pi} F(\t) e^{-ik\t} d\t$$ are the Fourier coefficients of $F$. Let $C_\ell (F)$ be the $\ell$-th cumulant of $S_N(F)=\sum_{j=1}^N F(\t_j)$. Then we have
		\begin{align}
			C_1(F) & =\hat{F}(0)\cdot N, \\
			\left|	C_2(F)-\sum_{k\in\mathbb Z} |k||\hat F(k)|^2  \right| & \le \sum_{|k|>N/2} |k||\hat F(k)|^2 , \label{eq sosh c2}
		\end{align}
		and for $\ell \ge 3$
		\begin{align}
			\left|	C_\ell(F) \right| \ll_\ell \sum_{\substack{k_1+\cdots +k_\ell=0\\ |k_1|+\cdots+|k_\ell|>N}} |k_1||\hat F(k_1)\cdots \hat F(k_\ell)|. \label{eq sosh c ell}
		\end{align}
	\end{prop} 
	
	We shall apply Proposition~\ref{prop sosh} to 
	\begin{align*}
		F(\t)=F_{N, L, u, v} (\t) = u g(\t) + v h(\t)
	\end{align*}
	for every $(u,v)\in \mathbb R^2$, and thus,
	\begin{align*}
		S_N(F) = 	uS_N(g)+vS_N(h).
	\end{align*}
	Since the normal distribution is determined by cumulants, to prove~\eqref{eq univariate} it suffices to prove that 
	\begin{align}\label{eq c to prove}
		C_1(F) \to 0, \quad C_2(F)\to \frac{u^2+v^2}{8}, \quad \text{ and } \quad C_\ell(F)\to 0 \text{ for } \ell \ge 3
	\end{align}
	as $N\to \infty$.
	
	We start by computing Fourier coefficients of $f$. Recall that
	\begin{align*}
		f(\t) =\frac{L}{N} \cdot \frac{1}{(1-L/N)-e^{i\t}}.
	\end{align*}
	It follows easily that
	\begin{align*}
		\hat{f}(k) = \left\{
		\begin{array}{ll}
			0, &\text{ if } k\ge 0, \\[1ex]
			\frac{-L}{N} \cdot \left(1-\frac{L}{N}\right)^{-k-1}, & \text{ if } k< 0.
		\end{array}
		\right.
	\end{align*}
	Note that here all $\hat{f}(k)$ are real. From this we deduce the Fourier coefficents for $g$ and $h$:
	\begin{align*}
		\hat{g}(k)= \frac{1}{2} \left( \hat{f}(k)+ \overline{\hat{f}(-k) } \right)= \left\{
		\begin{array}{ll}
			0, &\text{ if } k= 0, \\[1ex]
			\frac{-L}{2N} \cdot \left(1-\frac{L}{N}\right)^{|k|-1}, & \text{ if } k \ne 0,
		\end{array}
		\right.
	\end{align*}
	and
	\begin{align*}
		\hat{h}(k)= \frac{1}{2i} \left( \hat{f}(k)- \overline{\hat{f}(-k) } \right)= \left\{
		\begin{array}{ll}
			0, &\text{ if } k= 0, \\[1ex]
			\frac{-L}{2iN} \cdot \left(1-\frac{L}{N}\right)^{|k|-1}, & \text{ if } k< 0,  \\[1ex]
			\frac{L}{2iN} \cdot \left(1-\frac{L}{N}\right)^{|k|-1}, & \text{ if } k > 0.
		\end{array}
		\right.
	\end{align*}
	Since $F=ug+vh$, we have
	\begin{align}\label{eq F hat}
		\hat{F}(k) =  u	\hat{g}(k) + v	\hat{h}(k) =  \left\{
		\begin{array}{ll}
			0, &\text{ if } k= 0, \\[1ex]
			\frac{-L}{2N} \cdot \left(1-\frac{L}{N}\right)^{|k|-1}\cdot (u-iv), & \text{ if } k< 0,  \\[1ex]
			\frac{L}{2N} \cdot \left(1-\frac{L}{N}\right)^{|k|-1}\cdot (u+iv), & \text{ if } k > 0.
		\end{array}
		\right.
	\end{align}
	
	We first estimate the two sums in~\eqref{eq sosh c2}. A staightforward computation shows
	\begin{align*}
		\sum_{k\in\mathbb Z} |k||\hat F(k)|^2  = \frac{u^2+v^2}{2} \left(\frac{L}{N}\right)^2  \left(1-\frac{L}{N}\right)^{-2}  \sum_{k\ge 1} k \left(1-\frac{L}{N}\right)^{2k}.
	\end{align*}
	Denote temporarily $A(x)=  \sum_{k\ge 1} k x^{2k}$ and $B(x) =  \sum_{k\ge 1} x^{2k}$. For $0<x<1$ we have $B(x)=(1-x^2)^{-1} -1$. Differentiating $B(x)$ yields $A(x) = x^{2} (1-x^2)^{-2}$ for $0<x<1$. Therefore, letting $x=1-L/N$ in the above equation we obtain
	\begin{align*}
		\sum_{k\in\mathbb Z} |k||\hat F(k)|^2  = \frac{u^2+v^2}{2} \frac{1}{\left(2-\frac{L}{N}\right)^{2} },
	\end{align*}
	which is finite for fixed $N, L, u$ and $v$. Therefore, Proposition~\ref{prop sosh} applies to our function $F$.
	From the above equation we also see
	\begin{align} \label{eq c2 1}
		\lim_{N\to \infty} 	\sum_{k\in\mathbb Z} |k||\hat F(k)|^2  =  \frac{u^2+v^2}{8}. 
	\end{align}
	A similar treatment for the second sum in~\eqref{eq sosh c2} gives
	\begin{align*}
		\sum_{|k|>N/2} |k||\hat F(k)|^2 
		& = \frac{u^2+v^2}{2} \left(\frac{L}{N}\right)^2  \left(1-\frac{L}{N}\right)^{-2}  \sum_{k\ge N/2} k  \left(1-\frac{L}{N}\right)^{2k} \\
		& = \frac{u^2+v^2}{2} \left(\frac{L}{N}\right)^2  \left(1-\frac{L}{N}\right)^{-2}   \frac{\left(1-\frac{L}{N}\right)^{2M} }{\left(\frac{L}{N}\right)^{2}\left(2-\frac{L}{N}\right)^{2}  } \left(M+(1-M)\left(1-\frac{L}{N}\right)^{2} \right),
	\end{align*}
	where $M$ is the least integer greater than $N/2$. Since we assume $L=o(N)$, there is no harm to assume $L/N<1/2$, say. Thus, it is not difficult to see that the above is
	\begin{align*}
		&	\ll_{u,v} \left(1-\frac{L}{N}\right)^{2M+2} + \left(1-\frac{L}{N}\right)^{2M} \cdot M \cdot \frac{L}{N}\\
		&	\ll_{u,v} \left(1-\frac{L}{N}\right)^{2M} \cdot L \\
		&	\ll_{u,v} \left(1-\frac{L}{N}\right)^{N} \cdot L \\
		& 	\ll_{u,v} e^{-L} L.
	\end{align*}
	Since $L\to \infty$ with $N$, we have
	\begin{align} \label{eq c2 2}
		\lim_{N\to \infty}		\sum_{|k|>N/2} |k||\hat F(k)|^2  =0.
	\end{align}
	Combining \eqref{eq c2 1} and \eqref{eq c2 2} we conclude that
	\begin{align}
		\label{eq c2}
		C_2(F) \to \frac{u^2+v^2}{8} 
	\end{align}
	as $N\to \infty$.
	
	Moreover, from \eqref{eq F hat} and Proposition~\ref{prop sosh} it follows immediately that
	\begin{align}\label{eq c1}
		C_1(F) =0.
	\end{align}
	In view of \eqref{eq c to prove}, \eqref{eq c1}, \eqref{eq c2} and \eqref{eq sosh c ell}, it only remains to prove 
	\begin{align*}
		\sum_{\substack{k_1+\cdots +k_\ell=0\\ |k_1|+\cdots+|k_\ell|>N}} |k_1||\hat F(k_1)\cdots \hat F(k_\ell)| \to 0
	\end{align*}
	as $N\to \infty$, for each $\ell \ge 3$. From~\eqref{eq F hat} we have
	\begin{align*}
		|\hat{F}(k)|  \left\{
		\begin{array}{ll}
			=	0, &\text{ if } k= 0, \\[1ex]
			\ll_{u,v}	\frac{L}{N} \cdot \left(1-\frac{L}{N}\right)^{|k|}, & \text{ if } k \ne 0.
		\end{array}
		\right.
	\end{align*}
	Observe that 
	\begin{align*}
		\sum_{\substack{k_1+\cdots +k_\ell=0\\ |k_1|+\cdots+|k_\ell|>N}} |k_1||\hat F(k_1)\cdots \hat F(k_\ell)|
		 \ll_{\ell}   \sum_{\substack{k_1+\cdots +k_\ell=0\\ |k_1|+\cdots+|k_\ell|>N \\ |k_1| \ge |k_2| , \dots, |k_\ell|}} |k_1||\hat F(k_1)\cdots \hat F(k_\ell)| ,
	\end{align*}
and the conditions in the last sum  imply that $ |k_1| > N/\ell $. Thus, the above sum is
	\begin{align*}
		& \ll_{\ell} \sum_{|k| > N/\ell} |k| |\hat{F}(k)| \sum_{\substack{|k_2|\le |k|, ..., |k_\ell|\le |k| \\ k+k_2+\cdots +k_\ell =0 } } |\hat{F}(k_2)\cdots \hat{F}(k_\ell)|.
	\end{align*}	 
	Plug in the bounds for $|\hat{F}(k)|$, and note that the inner sum condition $k+k_2+\cdots +k_\ell =0 $ implies $ |k_2|+ \cdots + |k_\ell|  \ge |k| $. Thus, the above is
	\begin{align}	 \label{eq bound c ell}
		\notag	 & \ll_{\ell, u, v} \left(\frac{L}{N}\right)^\ell \cdot \sum_{|k| > N/\ell} |k|  \left(1-\frac{L}{N}\right)^{|k|}  \sum_{\substack{|k_2|\le |k|, ..., |k_\ell|\le |k| \\ k+k_2+\cdots +k_\ell =0 } } \left(1-\frac{L}{N}\right)^{|k_2|+ \cdots + |k_\ell|} 
		\\
		\notag	  	 & \ll_{\ell, u, v} \left(\frac{L}{N}\right)^\ell \cdot \sum_{|k| > N/\ell} |k|  \left(1-\frac{L}{N}\right)^{|k|}  \sum_{\substack{|k_2|\le |k|, ..., |k_{\ell-1}|\le |k| } }  \left(1-\frac{L}{N}\right)^{|k|}
		\\
		\notag & \ll_{\ell, u, v} \left(\frac{L}{N}\right)^\ell \cdot \sum_{|k| > N/\ell} |k|  \left(1-\frac{L}{N}\right)^{2|k|}  \sum_{\substack{|k_2|\le |k|, ..., |k_{\ell-1}|\le |k| } }  1
		\\
		\notag & \ll_{\ell, u, v} \left(\frac{L}{N}\right)^\ell \cdot \sum_{|k| > N/\ell} |k|^{\ell-1}  \left(1-\frac{L}{N}\right)^{2|k|}  
		\\	 
		&  \ll_{\ell, u, v} \left(\frac{L}{N}\right)^\ell \cdot \sum_{k> N/\ell} k^{\ell-1}  \left(1-\frac{L}{N}\right)^{2k}  .
	\end{align}
	Let $D_n(y) = \sum_{k\ge n} y^k= y^n(1-y)^{-1}$ for $0<y<1$. Differentiating $\ell-1$ times with respect to $y$, we have, for fixed $\ell$ and for $n>2\ell$, that
	\begin{align*}
		\sum_{k\ge n} k^{\ell-1} y^k & \le 	\sum_{k\ge n} k^{\ell-1} y^{k-\ell+1} \\
		& \ll_{\ell} \left(\frac{d}{dy}\right)^{\ell-1} D_n(y) \\
		& = \sum_{j=0}^{\ell-1} \binom{\ell-1}{j} \cdot \left(\frac{d}{dy}\right)^{j} y^n \cdot  \left(\frac{d}{dy}\right)^{\ell-1-j}\frac{1}{1-y}\\
		& \ll_{\ell}\sum_{j=0}^{\ell-1} n^j y^{n-j}\frac{1}{(1-y)^{\ell-j}}.
	\end{align*}
	Plug in $y=(1-L/N)^2$ and $n=$ the least integer $>N/\ell$, and let $N$ be sufficiently large. We obtain
	\begin{align*}
		\sum_{k> N/\ell} k^{\ell-1}  \left(1-\frac{L}{N}\right)^{2k} 
		& \ll_{\ell} \sum_{j=0}^{\ell-1} \left(\frac{N}{\ell}\right)^j \left(1-\frac{L}{N}\right)^{2(n-j)} \frac{1}{\left(2-\frac{L}{N}\right)^{\ell-j} \left(\frac{L}{N}\right)^{\ell-j}}\\
		& \ll_{\ell} \sum_{j=0}^{\ell-1} N^j \left(1-\frac{L}{N}\right)^{2n} \left(\frac{L}{N}\right)^{j-\ell}\\	
		& \ll_{\ell} N^\ell \left(1-\frac{L}{N}\right)^{2N/\ell} L^{-1}.
	\end{align*}
	From this and \eqref{eq bound c ell} it follows that
	\begin{align*}
		\sum_{\substack{k_1+\cdots +k_\ell=0\\ |k_1|+\cdots+|k_\ell|>N}} |k_1||\hat F(k_1)\cdots \hat F(k_\ell)|
		& \ll_{\ell,u,v}  \left(1-\frac{L}{N}\right)^{2N/\ell} L^{\ell-1} \\
		& \ll_{\ell,u,v}  e^{-2L/\ell} L^{\ell-1} 
	\end{align*}
	which tends to $0$ as $L\to \infty$ (or $N\to \infty$). This finishes the proof of Theorem~\ref{thm clt}. 
	
	\section{acknowledgments}
	
	I would like to thank Joseph Najnudel and  Ashkan Nikeghbali for a clarification of a result in~\cite{chhaibi2017ratios}. I also thank the referee for a very careful reading of the paper and for making many helpful suggestions.

	

	
\end{document}